\documentclass[12pt,a4paper]{article}
\usepackage[width=19cm,height=23cm]{geometry}
\usepackage{amsmath,amssymb,amsfonts,amsthm,enumerate,float}
\usepackage{mathrsfs,amsbsy,bm,mathtools}
\usepackage{indent first}
\usepackage{scalerel}
\DeclareMathOperator*{\Bigcdot}{\scalerel*{\cdot}{\bigodot}}

\newtheorem{Thm}{Theorem}[section]%
\newtheorem{Lem}[Thm]{Lemma}%
\newtheorem{Quest}[Thm]{Question}%

\theoremstyle{definition} 

\title{\bf Universal 
Entire Curves in
\\
 Projective Spaces with Slow Growth}
\author{Zhangchi Chen, Dinh Tuan Huynh, Song-Yan Xie}

\begin{document}
\maketitle

\abstract{
We construct explicit universal entire curves $h: \mathbb{C}\rightarrow \mathbb{C}\mathbb{P}^n$ whose Nevanlinna characteristic functions grow slower than any preassigned transcendental growth rate. Moreover, we can make $h$ to be hypercyclic for translation operations along any given countable  directions.
}



\section{\bf Introduction}
About one century ago,
in the space of entire functions, Birkhoff \cite{Birkhoff-1929} constructed some {\em universal} ones $h:\mathbb{C}\rightarrow\mathbb{C}$, which can approximate {\em any} entire function $g$ by translations of $h$. Precisely, given any compact set $K$ in $\mathbb{C}$, there is a sequence $(a_i)_{i\geqslant 1}$ of complex numbers such that $h(\Bigcdot+a_i)$
converges to $g(\Bigcdot)$ uniformly on $K$.

Such universality also appears  in the space  $\mathscr{M}_n$ of entire curves $f\colon\mathbb{C}\rightarrow\mathbb{C}\mathbb{P}^n$ in projective $n$-spaces, endowed with the topology of uniform convergence on compact sets.  Let $\mathscr{U}_n$ denote the subset of $\mathscr{M}_n$ consisting of
{\em universal entire curves}  defined likewise.

In this paper, we answer a question of Dinh-Sibony~\cite[Problem~9.1]{Dinh-Sibony-2020} in an extended version.

\begin{Quest}
	\label{Q 1.1}
Find  minimal growth rate of the Nevanlinna characteristic function $T_h(r)$ for $h\in\mathscr{U}_n$, where
\begin{equation}\label{eqn-char}
\normalfont T_h(r):=\int_{1}^{r}\dfrac{\text{d}t}{t}\int_{|z|<t}h^*{\omega_{{\sf FS}}}.
\end{equation}
\end{Quest}
In the formula $\omega_{{\sf FS}}$ is the Fubini-Study form on  $\mathbb{C}\mathbb{P}^n$. The Nevanlinna characteristic function
$T_f(r)$
measures the complexity of $f$. For instance,
 $f$ is a rational function if and only if $T_f(r) \approx O(1)\cdot\log~r$ as $r\rightarrow +\infty$ (cf. e.g.~\cite[Theorem~2.5.28]{Noguchi-Winkelmann-2014}).

One partial result in this direction was obtained in \cite{Luh-Martirosian-2000}, in which some meromorphic function
$f$ was constructed for approximating any entire functions by its translations, having Nevanlinna characteristic function
$T_f(r)\approx q(r)\cdot(\log~r)^2$ for some preassigned positive continuous function $q(r)$ tending to infinity as $r\rightarrow +\infty$.

Any universal entire curve $f$ must be transcendental, which guarantees that~\cite{Noguchi-Winkelmann-2014, Ru-2021}
\[
\lim_{r\rightarrow +\infty}
\,
T_f(r)/\log r
=
+\infty.
\]
In this paper we provide the following optimal answer to Dinh-Sibony's Question~\ref{Q 1.1}. 
\begin{Thm}\label{thm-main}
For any $n\geqslant 1$, given any positive continuous function $\phi:\mathbb{R}_+\rightarrow\mathbb{R}_+$ tending to infinity $\lim\limits_{r\rightarrow+\infty}\phi(r)=+\infty$, there exists some universal entire curve $h$ with slow growth 
\begin{align}\label{thr}
T_h(r)\leqslant\phi(r)\cdot\log~r, \quad\forall\, r\geqslant 1.
\end{align}
\end{Thm}

In history,
Birkhoff's discovery of universal entire functions foreshadowed a new branch of functional analysis  called hypercyclic operators theory~\cite{Bayart-Matheron-2009, Grosse-Erdmann-Manguillot-2011}. Now we briefly introduce the key notions. 

A continuous linear operator $\mathsf{T}: X\rightarrow X$ acting on a topological vector space $X$ (over $\mathbb{C}$ or $\mathbb{R}$) is called {\em hypercyclic} if there is some element $x\in X$ whose orbit $\{\mathsf{T}^n(x)\}_{n\geqslant 1}$ is dense in $X$. Such a vector $x$ is called {\em hypercyclic} for $\mathsf{T}$. 

From today's hindsight, Birkhoff's result can be rephrased as follows: for any nonzero complex number $a$, the translation operator $\mathsf{T}_a(f)(z):=f(z+a)$ is hypercyclic on the space $\mathscr{H}(\mathbb{C})$ of entire functions, endowed with the topology of uniform convergence on compact sets.

In this paper  we  show the existence of universal entire curves $f$ with slow growth, which are moreover hypercyclic simultaneously for all nontrivial translations along given countable directions.
One key ingredient is a reminiscent of~\cite[Theorem~8]{Costakis-Sambarino-2004}.

\begin{Thm}\label{thm-main2} For any $n\geqslant1$, given any positive continuous function $\phi:\mathbb{R}_+\rightarrow\mathbb{R}_+$ tending to infinity, and given any countable set $E\subset [0,2\pi)$, there exists some universal entire curve $h$ satisfying
\begin{itemize}
\item small growth rate $T_h(r)\leqslant\phi(r)\cdot\log~r$, for all $ r\geqslant 1$;
\item  $h$ is hypercyclic for $\mathsf{T}_a$ for any nonzero complex number $a$ with argument in $E$.
\end{itemize}
\end{Thm}

Our proof is based on insight of Oka manifolds theory~(cf.~\cite{Oka-book}) and Nevanlinna theory. Precisely, in our setting $\mathbb{CP}^n$ is a special Oka manifold with a ``large'' open subset $\mathbb{C}^n$  having complex vector space structure. The idea is, firstly, to construct some countable meromorphic discs
$\{g^{[k]}: \mathbb{D}_{R_k} \rightarrow \mathbb{C}^n\}_{k\geqslant 1}$, each being a rational map having designed approximation property and decaying to zero near infinity $\lim_{z\rightarrow \infty}g^{[k]}(z)=\vec{0}\in \mathbb{C}^n$. This step  can be accomplished by using  Runge's approximation theorem,
or by Lemma~\ref{thm-runge}.
Secondly, we ``patch'' the discs $\{g^{[k]}\}_{k\geqslant 1}$ together by sophisticated translations, repetitions, and infinite summation.

\medskip\noindent
{\bf Convention.} We denote by
$\mathbb{D}_r$  the disc centered at the origin with radius $r>0$, and by $\mathbb{D}(a,r)$ the disc  centered at $a\in \mathbb{C}$ with radius $r$.

\medskip

\noindent
{\bf Acknowledgement:}
Xie is
partially supported by 
National Key R\&D Program of China Grant
No.~2021YFA1003100 and  NSFC Grant No.~12288201. Chen is supported in part by the Labex CEMPI (ANR-11-LABX-0007-01), the project QuaSiDy (ANR-21-CE40-0016), and China Postdoctoral Science Foundation (2023M733690). Huynh is funded by University of
Education, Hue University under grant number NCM. T.22 – 02. We thank Bin Guo (UCAS) for pointing out the reference \cite{Costakis-Sambarino-2004}. We are grateful to the referee for nice suggestions.

\section{Preparations}
\subsection{Nevanlinna Theory}

Let $[Z_0:\dots:Z_n]$ be homogeneous coordinates of $\mathbb{C}\mathbb{P}^n$. Denote by $H_0:=\{Z_0=0\}$ the first coordinate hyperplane. Let  $h:\mathbb{C}\rightarrow\mathbb{C}\mathbb{P}^n$  be  an entire curve not contained in $H_0$,
with reduced representation 
$h=[h_0:h_1:\dots:h_n]$. Denote by $n_h(r,H_0)$ the number of zeros of $h_0$ in $\{|z|\leqslant r\}$, counting multiplicities. The {\em counting function} $N_h(r,H_0)$ of $h$ with respect to $H_0$ is
\[
N_h(r,H_0):=\int_{t=1}^r\frac{n_h(t,H_0)}{t}\,
\text{d}t,
\]
and the {\em proximity function} $m_h(r,H_0)$ is
\[
m_h(r,H_0):=\int_{\theta=0}^{2\pi}\log\,\sqrt{1+|h_1/h_0|^2(re^{i\theta})+\dots+|h_n/h_0|^2(re^{i\theta})}\,\frac{\text{d}\theta}{2\pi}.
\]
By Nevanlinna theory (cf.~\cite{Ru-2021}),   Shimizu-Ahlfors' version~\eqref{eqn-char} of Nevanlinna's characteristic function can be rewritten as~\cite[Theorem 2.3.31]{Noguchi-Winkelmann-2014}
\begin{align}\label{fmt}
	T_h(r)
	=
m_h(r,H_0)+N_h(r,H_0)-m_h(1,H_0).
\end{align}

\subsection{Runge-Type Approximations}
Runge's approximation theorem (1885) states that any holomorphic function defined on a neighborhood of a compact set $K\subset\mathbb{C}$ can be uniformly approximated by  rational functions with poles outside $K$. 
For our purpose, we need the following analogue.

\begin{Lem}\label{thm-runge} For any $f\in\mathscr{M}_n$, for any error bound $\epsilon>0$, for any $N>0$, there exists some rational map $$\gamma=[p_0:p_1:\dots:p_n]\in\mathscr{M}_n,$$
	where all $p_i\in \mathbb{C}[z]$ for $0\leqslant i\leqslant n$,
	 such that
\begin{enumerate}
\item[(i)] the Fubini-Study distance $d_{\sf FS}\big(\gamma(z),f(z)\big)<\epsilon$ on $\overline{\mathbb{D}_N}$;
\item[(ii)] for each $j\neq0$,
$\deg(p_0)>\deg(p_j)$, whence \[\left|p_j(z)/p_0(z)\right|=O\left(1/|z|\right),\quad  z\rightarrow \infty.
\]
\end{enumerate}
\end{Lem}
\proof Take a reduced representation $f=[f_0:\dots:f_n]$, where each $f_j$ is holomorphic on $\mathbb{C}$. By rescaling, we can assume that $|f_0|^2+\dots+|f_n|^2\geqslant1$ on $\overline{\mathbb{D}_N}$.

Geometrically, the Fubini-Study distance  between two points $[a_0:\dots:a_n]$ and $[b_0:\dots:b_n]$ is, up to multiplying by some positive constant, the angle between two complex vectors $\vec{a}=(a_0,\dots,a_n)$ and $\vec{b}=(b_0,\dots,b_n)\in\mathbb{C}^{n+1}$. For $\vec{a}$ 
 away  from the origin, say $||\vec{a}||\geqslant 1$, and for
 $\vec{b}$ sufficiently near to $\vec{a}$, say
 $||\vec{b}-\vec{a}||<\mu\ll 1$, 
 we have
\[
d_{\sf FS}([a_0:\dots:a_n],[b_0:\dots:b_n])<\epsilon.
\]

Therefore, we approximate each $f_j$ by some polynomial $p_j$ such that $|p_j(z)-f_j(z)|<\frac{\mu}{2\sqrt{n+1}}$ uniformly on $\overline{\mathbb{D}_N}$. Whence the requirement (i) is satisfied.

For the requirement (ii), we can start with a nonzero polynomial $p_0(z)$, and then replace it 
by $p_0(z)\frac{(z+M)^Q}{M^Q}$ for some large $M, Q\gg 1$, so  that $\left|p_0(z)\frac{(z+M)^Q}{M^Q}-p_0(z)\right|<\frac{\mu}{2\sqrt{n+1}}$ on $\overline{\mathbb{D}_N}$. \qed

\section{Construction}
We now prove Theorem~\ref{thm-main2} by an explicit construction.

 Let $E=\{\theta_v\}_{v\geqslant 1}\subset[0,2\pi)$ be a countable set of angles. Let $\phi: \mathbb{R}_+\rightarrow \mathbb{R}_+$ be a positive continuous function tending to infinity
 $\lim\limits_{r\rightarrow +\infty} \phi(r)=+\infty$. After replacing $\phi(r)$ by $\hat{\phi}(r):=\min_{t\geqslant r}\phi(t)$, we can assume that $\phi$ is nondecreasing. 

\subsection{Preparing Model Curves Over Discs}
We select a
countable and dense subfield of $\mathbb{C}$, say $\mathbb{Q}_c:=\{z=x+\sqrt{-1}y\, :~x,y\in\mathbb{Q}\}$.   
We will make use of the following countable subset of rational curves
\begin{equation}
\label{R}
\mathscr{R}:=\{[p_0(z):p_1(z):\dots:p_n(z)]\,:
\,
p_0, p_1, \dots, p_n \in \mathbb{Q}_c[z];\,
  \deg(p_0)>\deg(p_j),\forall\, j\neq 0\},
\end{equation}
which can be enumerated as $$\left(\gamma^{[k]}=[p_0^{[k]}(z):p_1^{[k]}(z):\dots:p_n^{[k]}(z)]\right)_{k\geqslant 1}.$$ By Lemma~\ref{thm-runge}, noting that $\mathbb{Q}_c\subset \mathbb{C}$ is dense,
we can find a sequence of elements in $\mathscr{R}$ to approximate any given entire curve uniformly on compact sets $\overline{\mathbb{D}_N}$.

Since $p_0^{[k]}\not\equiv0$, we can rewrite
\[
\gamma^{[k]}(z)=[1:\gamma_1^{[k]}(z):\dots:\gamma_n^{[k]}(z)], \quad \gamma_j^{[k]}(z):=p_j^{[k]}(z)/p_0^{[k]}(z).
\]
Note that $\left|\gamma_j^{[k]}(z)\right|=O(1/|z|)$ as $z\rightarrow \infty$, since   $\deg(p^{[k]}_0)>\deg(p^{[k]}_j)$. Choose a large radius $\eta_k>0$ such that $\left|\gamma_j^{[k]}(z)\right|\leqslant 2^{-k}$ for all $|z|>\eta_k$ and for all $j=1,\dots,n$.

A natural idea is constructing $h=[1:h_1:\dots:h_n]$ where $h_j=\sum_{k\geqslant 1} \gamma^{[k]}_j(z-a_k)$ for some sequence $(a_k)_{k\geqslant 1}$ of complex numbers growing sufficiently fast $$1\ll|a_1|\ll|a_2|\ll |a_3|\ll \cdots,$$ so that
\begin{itemize}
\item on each disc $\overline{\mathbb{D}(a_k,\eta_k)}$, the error term $\sum_{\ell\neq k}\gamma^{[\ell]}_j(z-a_{\ell})$ is negligible
\[
\Big|\sum_{\ell\neq k}\gamma^{[\ell]}_j(z-a_{\ell})
\Big|
\leqslant 
2^{-k+1},\quad
\forall\, j=1,\dots, n,\quad\forall\,
k\geqslant 1;
\]
\item outside all these discs, the sum $h_j$ has norm bounded by $1$.
\end{itemize}

However, this $h$ is not obviously universal due to the presence of error terms. To prove universality, we make the following manipulations.
\medskip

The Cartesian product
\[
\mathscr{A}:=\mathscr{R}\times\mathbb{Z}_+\times E=\left\{\big((\gamma^{[k]},u,\theta_v\big);~k,u,v\in\mathbb{Z}_+ \right\}
\]
is countable. We enumerate $\mathscr{A}$, for instance by lexicographic order of the index $(k,u,v)\in\mathbb{Z}_+^3$, as
\[
\mathscr{A}=(A_k)_{k\geqslant 1}, \quad A_k:=(g^{[k]},w_k,\rho_k).
\]
For any $\gamma\in\mathscr{R}$, for any $\theta\in E$, there is an infinite subsequence $A_{k_m}=(g^{[k_m]},w_{k_m},\rho_{k_m})$ for $m\geqslant 1$ such that $g^{[k_m]}=\gamma$, $\rho_{k_m}=\theta$, and $w_{k_m}$ tends to infinity as $m\rightarrow+\infty$.

Let
\[
h(z):=[1:h_1(z):\dots:h_n(z)],  \quad h_j(z):=\sum\limits_{k=1}^{+\infty}g^{[k]}_j(z-a_k).
\]
for some $a_k\in\mathbb{C}$ to be chosen subsequently by a sophisticated algorithm~\eqref{algorithm}. The following lemmas show that if the sequence $(|a_k|)_{k\geqslant 1}$ grows fast enough, we have the following properties:
\begin{enumerate}
\item[(1)]  (convergence) $h_j$ is well-defined meromorphic function for each $j=1,\dots,n$; 
\item[(2)]  (universality) $h$ is a universal; 
\item[(3)] (growth rate) $T_h(r)\leqslant \phi(r)\cdot\log~r$, for all $ r\geqslant 1$.
\end{enumerate}

\subsection{Proof of Convergence and Universality}

Since $\left|g_j^{[k]}(z)\right|=O(1/|z|)$ for all $1\leqslant j\leqslant n$, there are some $\delta_k>0$ and $C_k>0$ such that $\left|g_j^{[k]}(z)\right|<\frac{C_k}{|z|}$ for any $|z|>\delta_k$. Set $R_k:=\delta_k+2^{k}\,C_k$. Hence
\begin{equation}
	\label{ineq 1}
\left|g_j^{[k]}(z)\right|
<\frac{C_k}{|z|}<2^{-k}\cdot \frac{R_k}{|z|}<2^{-k}, \quad \forall\,|z|>R_k,\quad \forall\, 1\leqslant j\leqslant n.
\end{equation}
We demand that all the discs $\mathbb{D}(a_k, R_k)$ for $k\geqslant 1$ are sufficiently far way from each other, in order to have the following good estimates.

\begin{Lem}\label{lem-sep} If
\begin{align}\label{condition-sep}
|a_k|-|a_{k-1}|-R_k-R_{k-1}>(R_1+\dots+R_{k-1})(k-1)2^k, \quad\forall\,k\geqslant 2,
\end{align}
then for any $j\in\{1,\dots,n\}$:
\begin{enumerate}
\item outside $\bigcup\limits_{\ell}\mathbb{D}(a_{\ell},R_{\ell})$, the function $h_j(z)=\sum\limits_{k=1}^{+\infty}g_j^{[k]}(z-a_{k})$ has bound $|h_j(z)|< 1$;
\item in each disc $\overline{\mathbb{D}(a_k,R_k)}$, the error term $\epsilon_j^{[k]}(z):=\sum\limits_{\ell\neq k}g_j^{[\ell]}(z-a_l)$ has bound $\left|\epsilon_j^{[k]}\right|< 2^{-k+1}$.
\end{enumerate}
\end{Lem}
\proof The first estimate follows by summing up~\eqref{ineq 1} for $k\geqslant 1$.
Now we prove the second estimate.

For each ${\ell}<k$, the distance between $\mathbb{D}(a_{\ell},R_{\ell})$ and $\overline{\mathbb{D}(a_k,R_k)}$ is at least $ R_{\ell}\,(k-1)\,2^{k}$ by~\eqref{condition-sep}. Hence by~\eqref{ineq 1} we have
\[
\left|g_j^{[{\ell}]}(z-a_{\ell})\right|<2^{-{\ell}}\frac{R_{\ell}}{R_{\ell}\,(k-1)\,2^{k}}<\frac{2^{-k}}{k-1}, \quad\forall\, z\in \mathbb{D}(a_k,R_k), \quad\forall\, {\ell}<k.
\]
Lastly, using the above estimate for $\ell<k$ and using~\eqref{ineq 1} for $\ell>k$, we receive
\[
\left|\epsilon_j^{[k]}(z)\right|
\leqslant
\sum\limits_{{\ell}<k}\left|g_j^{[{\ell}]}(z-a_{\ell})\right|+\sum\limits_{{\ell}>k}\left|g_j^{[{\ell}]}(z-a_{\ell})\right|
< \sum\limits_{{\ell}<k}\frac{2^{-k}}{k-1}+\sum\limits_{{\ell}>k}2^{-{\ell}}=2^{-k}+2^{-k}.
\]
\qed

Consequently, for each $j\in\{1,\dots,n\}$, the infinite sum $h_j(z)$ is a well-defined meromorphic function, and all its poles lie in $\bigcup\limits_{k} \mathbb{D}(a_k,R_k)$.

\begin{figure}[htbp]
\begin{center}
\includegraphics[width=0.6\linewidth]{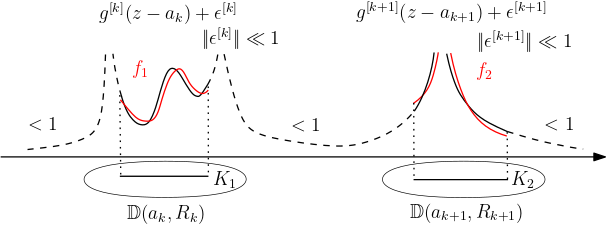} 
\caption{The dotted universal curve $h$, locally expressed by $g+\epsilon$, approximates a curve $f_1$ on $K_1$ and a curve $f_2$ (intersecting $H_0$) on $K_2$}
\label{fig:dn}
\end{center}
\end{figure}

\begin{Lem} If the condition~\eqref{condition-sep}  holds, then $h$ is universal.
\end{Lem}
\proof It suffices to prove that for any entire curve $f\in\mathscr{M}_n$, for any error bound $\epsilon>0$, and for any compact set $K$ in $\mathbb{C}$, there is some translation $\mathsf{T}_c$ given by large $c\in\mathbb{C}$, such that 
\[
d_{\sf FS}\big(f(z), \mathsf{T}_c(h)(z)\big)
=
d_{\sf FS}\big(f(z), h(z+c)\big)
<\epsilon,
\quad \forall\,z\in K.
\]

Take a large disc $\overline{\mathbb{D}_N}$
containing $K$.
 By Lemma~\ref{thm-runge},
 recalling~\eqref{R},
 we  can find some rational curve $\gamma=[p_0:p_1:\dots:p_n]\in\mathscr{R}$ such that
\[
d_{\sf FS}\left(f(z),\gamma(z)\right)<\frac{\epsilon}{2}, \quad \forall\, z\in \overline{\mathbb{D}_N}.
\]
By rescaling, we can moreover assume that $\sum\limits_{j=0}^n|p_j(z)|^2\geqslant 1$ on $\overline{\mathbb{D}_N}$.

By our construction, $\gamma$ repeats infinitely many times in $\left(g^{[k]}\right)_{k\geqslant 1}$.
 Take $k\gg 1$ such that
\begin{itemize}
\item the model curve $g^{[k]}=\gamma=[1:p_1/p_0:\dots:p_n/p_0]=:[1:g^{[k]}_1:\dots:g^{[k]}_n]$;
\item the error term $\left|\epsilon_j^{[k]}(z)\right|\leqslant\frac{\delta}{M}$ on $\overline{\mathbb{D}(a_k,R_k)}$ for all $j\in\{1,\dots,n\}$, 
where  $M:=\max_{z\in\overline{\mathbb{D}_N}}\,\{|p_0(z)|\}$, and
where $\delta>0$ is a sufficiently mall number such that
\[
d_{\sf FS}([\vec{a}],[\vec{b}])<\epsilon/2, \ \ \ \ \forall~ \vec{a},\vec{b}\in\mathbb{C}^{n+1}, |\!|\vec{a}|\!|\geqslant 1, |\!|\vec{b}-\vec{a}|\!|<\delta;
\]
\item the radius $R_k>N$.
\end{itemize}
Hence $\left|\epsilon_j^{[k]}(z+a_k)\right|\leqslant\frac{\delta}{M}$ on $\overline{\mathbb{D}_{R_k}}$. Rewrite
\[
\aligned
h(z+a_k)&=[1:h_1(z+a_k):\dots:h_n(z+a_k)]\\
&=[1:g^{[k]}_1(z)+\epsilon^{[k]}_1(z+a_k):\dots:g^{[k]}_n(z)+\epsilon^{[k]}_n(z+a_k)]\\
&=[p_0(z):p_1(z)+\epsilon^{[k]}_1(z+a_k)\,p_0(z):\dots:p_n(z)+\epsilon^{[k]}_n(z+a_k)\,p_0(z)],
\endaligned
\]
where $\left|\epsilon_j^{[k]}(z+a_k)\cdot p_0(z)\right|\leqslant \frac{\delta}{M}\cdot M= \delta$ on $\overline{\mathbb{D}_N}$. Since $\sum\limits_{j=0}^n|p_j(z)|^2\geqslant 1$, using the argument in the proof of  Lemma~\ref{thm-runge}, this sufficiently small $\delta$ ensures that
\[
d_{\sf FS}\left(h(z+a_k),\gamma(z)\right)<\frac{\epsilon}{2}, \quad \forall\,z\in\overline{\mathbb{D}_N}.
\]
Hence $d_{\sf FS}\left(h(z+a_k),f(z)\right)<\epsilon$ on $\overline{\mathbb{D}_N}$.\qed

\subsection{Estimate of the Growth Rate}

Now we prove that if the sequence $(|a_k|)_{k\geqslant 1}$ grows fast enough, \eqref{thr} holds.

Denote by $n_{k,j}$ the number of poles of the rational function $g_j^{[k]}$, counting with multiplicities. Let $n_k:=\sum_{j=1}^n n_{k,j}$.

We demand that $|a_k|>R_k+1$ for all $k\geqslant 1$.

\begin{Lem}\label{lem-growth} If
\begin{align}\label{condition-growth}
\left\{
\begin{aligned}
|a_k|-|a_{k-1}|-R_k-R_{k-1}&>(R_1+\dots+R_{k-1})(k-1)2^k, & \quad\forall\,k\geqslant 2,\\
\phi(|a_k|-R_k)&>(n_1+\dots+n_k+\sqrt{n+1})\,\tfrac{\log(|a_k|+R_k)}{\log(|a_k|-R_k)}, & \quad\forall\,k\geqslant 1,
\end{aligned}
\right. 
\end{align}
then $T_h(r)\leqslant \phi(r)\cdot\log~r$, for all $r\geqslant  \max\{|a_1|-R_1,e\}$.
\end{Lem}

\begin{figure}[htbp]
\begin{center}
\includegraphics[width=0.7\linewidth]{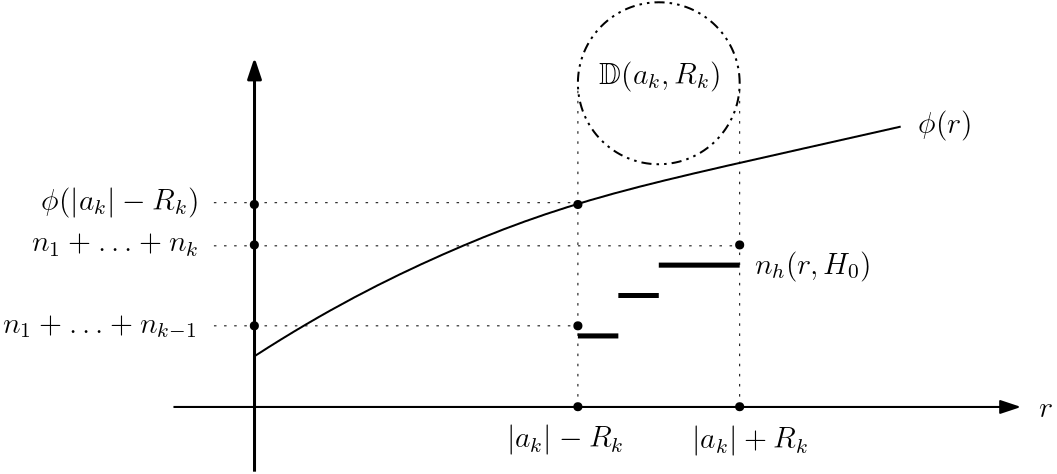} 
\caption{When $(|a_k|)_{k\geqslant 1}$ grows sufficiently fast, $n_h(r,H_0)$ is smaller than $\phi(r)$}
\label{fig:nhr}
\end{center}
\end{figure}

\proof By~\eqref{fmt}, for $r\geqslant 1$, we have
\begin{equation}\label{fmtt}
T_h(r)=N_h(r,H_0)+m_h(r,H_0)-m_h(1,H_0)\leqslant N_h(r,H_0)+m_h(r,H_0).
\end{equation}

\noindent{\bf Case 1:} There is some $k\in\mathbb{Z}_+$ such that $|a_k|+R_k<r<|a_{k+1}|-R_{k+1}$. In this case the circle $\{|z|=r\}$ does not intersect with any disc $\overline{\mathbb{D}(a_\ell,R_\ell)}$, $\ell\geqslant 1$.

By Lemma~\ref{lem-sep}, on this circle $|h_j(z)|\leqslant 1$ for all $j\in\{1,\dots,n\}$. Hence the proximity function
\begin{equation}\label{mhr}
m_h(r,H_0)\leqslant \sqrt{n+1}.
\end{equation}

Although the expression $h(z)=[1:h_1(z):\dots:h_n(z)]$ is not a reduced representation, nevertheless, we have the clear estimate $n_h(r,H_0)\leqslant n_1+\dots+n_k$. Therefore
\begin{equation}
	\label{Nhr}
N_h(r,H_0)\leqslant (n_1+\dots+n_k)\cdot \log~r.
\end{equation}
Summing up~\eqref{mhr} and~\eqref{Nhr}, and using~\eqref{fmtt}, we get
\begin{align}\label{ak+rk}
T_h(r)&\leqslant (n_1+\dots+n_k)\cdot \log~r+\sqrt{n+1}\\
\nonumber
\text{[use~\eqref{condition-growth} and $r\geqslant e$]}\quad&<\phi(r)\cdot\log~r.
\end{align}

\medskip

\noindent{\bf Case 2:} There is some $k\in\mathbb{Z}_+$ such that $\big|r-|a_k|\big|\leqslant R_k$. In this case, the circle $\{|z|=r\}$ intersects only one disc $\overline{\mathbb{D}(a_k,R_k)}$. 

By definition, $T_{h}(r)$ is nondecreasing. So
\[
\aligned
T_{h}(r)&\leqslant T_{h}(|a_k|+R_k)\\
\text{[use~\eqref{ak+rk}]}\quad&\leqslant(n_1+\dots+n_k+\sqrt{n+1})\cdot\log(|a_k|+R_k)\\
\text{[use~\eqref{condition-growth}]}\quad&<\phi(|a_k|-R_k)\cdot\log(|a_k|-R_k)\\
[\text{use~} |a_k|-R_k\leqslant r]\quad&\leqslant\phi(r)\cdot\log~r.
\endaligned
\]\qed

Now we estimate the growth rate for $r\in[1,|a_1|-R_1]$. The following Lemma tells that the Nevanlinna characteristic function $T_h(r)$ of an ``almost constant'' curve $h$ has very slow growth for small $r$.
\begin{Lem}\label{lem-control-T} Let $\phi:\mathbb{R}_+\rightarrow\mathbb{R}_+$ be a positive continuous nondecreasing function. For any $r_0>1$, there is $\epsilon_0>0$ such that if $h=[1:h_1(z):\dots:h_n(z)]$ for some meromorphic functions $h_j$ with $|h_j(z)|\leqslant\epsilon_0$ on $\overline{\mathbb{D}_{r_0+1}}$, then $T_h(r)\leqslant\phi(r)\cdot\log~r$, for all $r\in[1,r_0]$. 
\end{Lem}
\proof By Cauchy's integral formula
\[
h_j'(z)=\frac{1}{2\pi\,i}\int_{|\zeta|=r_0+1}\frac{h_j(\zeta)}{(\zeta-z)^2}\,\text{d}\zeta
, \quad \forall\,z\in\overline{\mathbb{D}_{r_0}},
\]
we receive the derivative estimate $|h_j'(z)|\leqslant (r_0+1)\,\epsilon_0$ on $\overline{\mathbb{D}_{r_0}}$. The Fubini-Study from on the affine chart $\mathbb{C}^n\subset\mathbb{C}\mathbb{P}^n$ with coordinates $(z_1,\dots,z_n)\mapsto[1:z_1:\dots:z_n]$, is
\[
\omega_{\sf FS}=\frac{i}{2\pi}\partial\bar{\partial}\log\left(1+\sum\limits_{j=1}^n|z_j|^2\right)=\frac{i}{2\pi}\left(\frac{\sum\limits_{j=1}^n\,\text{d}z_j\wedge\,\text{d}\overline{z}_j}{1+\sum\limits_{j=1}^n|z_j|^2}-\frac{\sum\limits_{j=1}^n\,\overline{z}_j\,\text{d}z_j}{1+\sum\limits_{j=1}^n|z_j|^2}\wedge\frac{\sum\limits_{j=1}^n\,z_j\,\text{d}\overline{z}_j}{1+\sum\limits_{j=1}^n|z_j|^2}\right).
\]
The pull-back of the Fubini-Study form is
\[
h^*\omega_{\sf FS}=\left(\frac{\sum\limits_{j=1}^n|h_j'|^2}{1+\sum\limits_{j=1}^n|h_j|^2}
-
\frac{\sum\limits_{j=1}^n\overline{h_j}h_j'}{1+\sum\limits_{j=1}^n|h_j|^2}\,\frac{\sum\limits_{j=1}^n\,h_j\overline{h_j}'}{1+\sum\limits_{j=1}^n|h_j|^2}
\right)\frac{i}{2\pi}\,\text{d}z\wedge \text{d}\overline{z}.
\]
Here the coefficient
\[
\left(\frac{\sum\limits_{j=1}^n|h_j'|^2}{1+\sum\limits_{j=1}^n|h_j|^2}
-
\frac{\sum\limits_{j=1}^n\overline{h_j}h_j'}{1+\sum\limits_{j=1}^n|h_j|^2}\,\frac{\sum\limits_{j=1}^n\,h_j\overline{h_j}'}{1+\sum\limits_{j=1}^n|h_j|^2}
\right)\leqslant\frac{\sum\limits_{j=1}^n|h_j'|^2}{1+\sum\limits_{j=1}^n|h_j|^2}\leqslant n(r_0+1)^2\,\epsilon_0^2
\]
on $\overline{\mathbb{D}_{r_0}}$. Hence for any $r\in[1,r_0]$, we have
\[
\aligned
T_h(r)\leqslant\int_{t=1}^{r}\frac{\text{d}t}{t}\int_{|z|\leqslant r_0}h^*\omega_{\sf FS}\leqslant \int_{t=1}^{r}\frac{\text{d}t}{t}\int_{|z|\leqslant r_0}n\,(r_0+1)^2\,\epsilon_0^2\frac{i}{2\pi}\,\text{d}z\wedge \text{d}\overline{z}=n\,r_0^2(r_0+1)^2\epsilon_0^2\,\log~r.
\endaligned
\]
We complete the proof by taking $\epsilon_0>0$ sufficiently small such that $n\,r_0^2(r_0+1)^2\epsilon_0^2\leqslant\phi(1)$. \qed

\begin{Lem}
	\label{dirty lemma} Take large $r_0>e$ such that $\phi(r)\cdot\log~r>\sqrt{n+1}$, for all $r>r_0$. Set $$\epsilon_0=\min\left\{1,\sqrt{\tfrac{\phi(1)}{n\,r_0^2\,(r_0+1)^2}}\right\}.$$ Under the condition (\ref{condition-growth}) and, if moreover
\begin{align}\label{condition-growth-3}
|a_k|>R_k/\epsilon_0+r_0+1, \quad \forall\,k\geqslant 1,
\end{align}
then $T_h(r)\leqslant \phi(r)\cdot\log~r$, for all $r\in[1,R_1-|a_1|)$.
\end{Lem}
\proof Recalling~\eqref{ineq 1}
\[
\Big|g^{[k]}_j(z)\Big|<2^{-k}\cdot \frac{R_k}{|z|}, \quad \forall\,|z|>R_k,
\] 
noting that
\[
|z-a_k|\geqslant R_k/\epsilon_0\geqslant R_k, \quad \forall\, |z|\leqslant r_0+1,
\]
we receive 
\[
\Big|g^{[k]}_j(z-a_k)\Big|\leqslant 2^{-k}\,\frac{R_k}{|z-a_k|}\leqslant2^{-k}\,\frac{R_k}{R_k/\epsilon_0}=2^{-k}\,\epsilon_0, \quad \forall\,|z|\leqslant r_0+1.
\]

Summing up the above inequality for $k\geqslant 1$, we get $|h_j(z)|\leqslant\epsilon_0$ on $\overline{\mathbb{D}_{r_0+1}}$. By Lemma~\ref{lem-control-T}, $T_h(r)\leqslant \phi(r)\cdot\log~r$, for $r\in[1,r_0]$.

For the remaining $r\in(r_0,\max\{|a_1|-R_1,e\})$, the proximity function $m_h(r,H_0)\leqslant\sqrt{n+1}$ since $|h_j(z)|\leqslant 1$ by Lemma~\ref{lem-sep}, and the counting function $N_h(r,H_0)=0$ by our construction.
Thus we conclude the proof by~\eqref{fmt}
\[
T_h(r)\leqslant N_h(r,H_0)+m_h(r,H_0)\leqslant \sqrt{n+1}<\phi(r)\cdot\log~r.\qedhere
\] 

\subsection{End of the Proof}
Now we summarize the first restrictions on the sequence $(a_k)_{k\geqslant 1}$:
\begin{equation}
	\label{algorithm}
\left\{
\begin{aligned}
|a_k|&>R_k/\epsilon_0+r_0+1, &\quad\forall\, k\geqslant 1,\\
|a_k|-|a_{k-1}|-R_k-R_{k-1}&>(R_1+\dots+R_{k-1})(k-1)2^k, &\quad\forall\, k\geqslant 2,\\
\phi(|a_k|-R_k)&>(n_1+\dots+n_k+\sqrt{n+1})\,\tfrac{\log(|a_k|+R_k)}{\log(|a_k|-R_k)}, &\quad\forall\,k\geqslant 1.
\end{aligned}
\right.
\end{equation}
In the first inequality, $\epsilon_0$ and $r_0$ are defined in Lemma~\ref{dirty lemma}.
Moreover,  we demand that each complex number
\[
a_k\in \mathbb{Z}_{+}\cdot e^{\sqrt{-1}\rho_k},
\quad
\forall\,k\geqslant 1,
\]
has the argument $\rho_k$. 

Clearly, such $a_k$ can be chosen
subsequently for $k=1, 2, 3, \dots$, 
since
\[
\lim\limits_{|a_k|\rightarrow+\infty}(n_1+\dots+n_k+\sqrt{n+1})\,\tfrac{\log(|a_k|+R_k)}{\log(|a_k|-R_k)}
<+\infty
=
\lim_{r\rightarrow +\infty} \phi(r).
\]
Thus the entire curve $h$ is universal with slow growth~\eqref{thr} in Theorem~\ref{thm-main} by the previous arguments.

 Lastly,  $h$ is hypercyclic for any translation operator $\mathsf{T}_{e^{\sqrt{-1}\theta}}$ with $\theta\in E$,
 since $\theta$ repeats infinitely many time in $\{\rho_k\}_{k\geqslant 1}$. Using the argument of a result \cite[Theorem~8]{Costakis-Sambarino-2004}, we conclude that $h$ is simultaneously hypercyclic for  any $\mathsf{T}_{a}$ where $a\in \mathbb{R}_+\cdot e^{\sqrt{-1}\theta}$ has argument $\theta$.

\setlength\parindent{0em}

{\scriptsize Zhangchi Chen, Academy of Mathematics and System Science, Morningside Center of Mathematics, Chinese Academy of Science, Beijing 100190, China}\\
{\bf\scriptsize zhangchi.chen@amss.ac.cn}

\medskip

{\scriptsize Dinh Tuan Huynh, Department of Mathematics, University of Education, Hue University, 34 Le Loi St., Hue City, Vietnam}\\
{\bf\scriptsize huynhdinhtuan@dhsphue.edu.vn}

\medskip

{\scriptsize Song-Yan Xie, Academy of Mathematics and System Science \& Hua Loo-Keng Key Laboratory
		of Mathematics, Chinese Academy of Sciences, Beijing 100190, China; 
		School of Mathematical Sciences, University of Chinese Academy of Sciences, Beijing
100049, China}\\
{\bf\scriptsize xiesongyan@amss.ac.cn}


\begin{thebibliography}{XL}{\scriptsize

{\bf\bibitem{Bayart-Matheron-2009}
{\rm Bayart}}, F.; {\rm Matheron}, É.:
{\em Dynamics of linear operators.}
Cambridge Tracts in Mathematics, 179. Cambridge University Press, Cambridge, 2009. xiv+337 pp. ISBN: 978-0-521-51496-5

\medskip

{\bf\bibitem{Birkhoff-1929}
{\rm Birkhoff}}, G.D.:
{\em D\'emonstration d'un th\'eor\`eme \'el\'ementaire sur les fonctions enti\`eres.}
C.R. Acad. Sci. Paris, 189:473--475, 1929.

\medskip

{\bf\bibitem{Dinh-Sibony-2020}
{\rm Dinh}} T-C.; {\rm Sibony} N.;
{\em Some open problems on holomorphic foliation theory.}
Acta Math. Vietnam., 45(1):103--112, 2020.

\medskip

{\bf\bibitem{Costakis-Sambarino-2004}
{\rm Costakis}} G.; {\rm Sambarino} M.:
{\em Genericity of wild holomorphic functions and common hypercyclic vectors.}
Adv. Math. 182 (2004), no. 2, 278--306.

\medskip

{\bf\bibitem{Oka-book}
	{\rm Forstneri\v{c}}}, Franc:
{\em Stein manifolds and holomorphic mappings.}
Springer, Cham; Ergebnisse der Mathematik und ihrer Grenzgebiete Vol. 56 (2nd edition). xiv+562 pp. 2017.


\medskip

{\bf\bibitem{Grosse-Erdmann-Manguillot-2011}
{\rm Grosse-Erdmann}} K-G.; {\rm Manguillot} A.P.:
{\em Linear chaos.}
Springer London; Universitext; 2011.

\medskip

{\bf\bibitem{Luh-Martirosian-2000}
{\rm Luh}}, W.; {\rm Martirosian}, V.:
{\em On the growth of universal meromorphic functions.}
Analysis (Munich) 20 (2000), no. 2, 137--147.

\medskip

{\bf\bibitem{Noguchi-Winkelmann-2014}
{\rm Noguchi}} J.; {\rm Winkelmann} J.:
{\em Nevanlinna theory in several complex variables and Diophantine approximation.}
Grundlehren der mathematischen Wissenschaften [Fundamental Principles of Mathematical Sciences], 350. Springer, Tokyo, (2014). xiv+416 pp. ISBN: 978-4-431-54570-5; 978-4-431-54571-2

\medskip

{\bf\bibitem{Ru-2021}
{\rm Ru}} M.:
{\em Nevanlinna theory and its relation to Diophantine approximation.}
Second edition [of MR1850002]. World Scientific Publishing Co. Pte. Ltd., Hackensack, NJ, (2021). xvi+426 pp. ISBN: 978-981-123-350-0; 978-981-123-351-7; 978-981-123-352-4
}
\end{thebibliography}
\end{document}